\documentclass[preprint,12pt,3p]{elsarticle}
\usepackage{amsmath}
\usepackage{amssymb}
\usepackage[dvipsnames]{xcolor}
\usepackage[utf8]{inputenc}
\usepackage{verbatim}
\usepackage[mathcal]{euscript}
\usepackage{graphicx}
\usepackage[english]{babel}
\usepackage{multicol}
\usepackage[utf8]{inputenc}  
\usepackage[T1]{fontenc}  
\usepackage{caption}
\usepackage[colorlinks=true ]{hyperref}
\usepackage[final]{pdfpages}
\usepackage{multicol}
\usepackage{float}
\usepackage{tikz}
\usepackage{tikz-3dplot} 
\usepackage{pdfpages}   
\usepackage{varioref} 
\usepackage{float}
\usepackage{lineno}
 \usepackage{multicol}
\usepackage{multirow}
\usetikzlibrary{quotes,arrows.meta}
\usetikzlibrary{shapes,arrows}
\numberwithin{equation}{section}  

\definecolor{camel}{rgb}{0.76, 0.6, 0.42}
\tdplotsetmaincoords{70}{120} 

\newtheorem{definition}{Definition}[section]

\newtheorem{Claim}[definition]{Claim}
\newtheorem{Remark}[definition]{Remark}
\newtheorem{example}[definition]{Example}

\newtheorem{proof}[definition]{Proof}
\newtheorem{Corollary}[definition]{Corollary}

\newcommand \bei {\begin{itemize}}
\newcommand \eei {\end{itemize}}
\newcommand \ubar u

\newcommand \del \partial

\newcommand \la \langle
\newcommand \ra \rangle 

\newcommand \auth    \textsc
\newcommand \be {\begin{equation}}
\newcommand \ee {\end{equation}}

\newcommand \bcor {\begin{Corollary}}
\newcommand \ecor {\end{Corollary}}

\newcommand \bpro {\begin{proof}}
\newcommand \epro {\end{proof}}
\newcommand \bdf {\begin{definition}}
\newcommand \edf {\end{definition}}
\newcommand \bex {\begin{example}}
\newcommand \eex {\end{example}}
\newcommand \bcl {\begin{Claim}}
\newcommand \ecl {\end{Claim}}
\newcommand \brm {\begin{Remark}}
\newcommand \erm {\end{Remark}}

\let\oldmarginpar\marginpar
\renewcommand\marginpar[1]{\-\oldmarginpar[\raggedleft\footnotesize #1]%
{\raggedright\footnotesize #1}}
\begin{document}
\include{TeXFigs}
 \begin{frontmatter}

\title{Implicit EXP-RBF techniques for modeling unsaturated flow through soils with water uptake by plant roots}

\author[mpu]{Mohamed Boujoudar}
\author[mpu,uottawa]{Abdelaziz Beljadid\corref{mycorrespondingauthor}}
\cortext[mycorrespondingauthor]{Corresponding author}
\ead{abdelaziz.beljadid@um6p.ma; abeljadi@uottawa.ca}
\author[LMA]{Ahmed Taik}
\address[mpu]{Mohammed~VI Polytechnic University, Morocco}
\address[uottawa]{ University of Ottawa, Canada}
\address[LMA]{University  Hassan~II, Morocco}
\begin{abstract}  
Modeling unsaturated flow through soils with water uptake by plan root has many applications in agriculture and water resources management. In this study, our aim is to develop efficient numerical techniques for solving the Richards equation with a sink term due to plant root water uptake.
The Feddes model is used for water absorption by plant roots, and the van-Genuchten model is employed for capillary pressure. 
We introduce a numerical approach that combines the localized exponential radial basis function (EXP-RBF) method for space and the second-order backward differentiation formula (BDF2) for temporal discretization. The localized RBF methods eliminate the need for mesh generation and avoid ill-conditioning problems. This approach yields a sparse matrix for the global system, optimizing memory usage and computational time. The proposed implicit EXP-RBF techniques have advantages in terms of accuracy and computational efficiency thanks to the use of BDF2 and the localized RBF method. Modified Picards iteration method for the mixed form of the Richards equation is employed to linearize the system. Various numerical experiments are conducted to validate the proposed numerical model of infiltration with plant root water absorption. The obtained results conclusively demonstrate the effectiveness of the proposed numerical model in accurately predicting soil moisture dynamics under water uptake by plant roots. The proposed numerical techniques can be incorporated in the numerical models where unsaturated flows and water uptake by plant roots are involved such as in hydrology, agriculture, and water management.
\end{abstract} 
\end{frontmatter}

\section{Introduction} 
%
%
Understanding the dynamics of water flows in unsaturated soils and their interactions with plant root systems are critical for agricultural efficiency, environmental sustainability, and climate predictions. These interactions play a crucial role in the hydrological cycle, particularly within the vadose zone. This zone has a significant impact on subsurface water resources and acts as a reservoir for soil moisture, which offers essential water for plant growth. The mechanism by which plants absorb water from the soil, known as root water uptake, has an important influence on the distribution of soil moisture. A comprehensive understanding of the plant root water uptake process is pivotal due its applications in agriculture, irrigation, and water management. Understanding this process can aid in developing strategies to enhance crop productivity, improve water use efficiency, and better manage limited water resources in various agricultural contexts.

The plant root water uptake is influenced by soil heterogeneity and root architecture \citep{FeddesR.A1978Sofw,feddes2004parameterizing}. Furthermore, environmental variables such as humidity and temperature affect the plant’s transpiration rate, which in turn impacts water uptake \cite{tardieu1998variability}. These complex interactions between soil properties, root characteristics, and environmental factors make root water uptake a challenging process to study and understand comprehensively. Robust coupled numerical models can be used as efficient tools to understand these complex interactions and study the impact of root water uptake on soil water distribution. These coupled models provide accurate simulations of water dynamics in the root zone and contribute to the enhancement of water resource management and informed decision-making in agriculture applications.


The infiltration of water into soils can be described at the continuum scale using Richards' equation~\citep{richards1931capillary}. This equation, derived from the combination of Buckingham-Darcy's law and the mass conservation equation \citep{bear1987theory}, governs the dynamics of water movement through unsaturated soils. This equation exhibits high non-linearity due to the relationship between hydraulic conductivity and capillary pressure which are dependent on the soil water content \citep{gardner1958some,brooks1964hydrau,van1980closed}. Water uptake by roots is commonly described using two main approaches: macroscopic and microscopic \citep{gardner1960dynamic,molz1968soil,molz1981models}. Microscopic models, which are physically based, describe water extraction at the individual root level and consider the radial flow of soil water towards specific roots \citep{gardner1960dynamic,hillel1975microscopic,hainsworth1986water}. These approaches require detailed information about root geometry. The second category employs an empirical macroscopic approach \citep{FeddesR.A1978Sofw,van1987numerical,albasha2015compensatory}, and considers the entire root system as a unified entity to account for the combined effects of individual roots. In the macroscopic approach, the plant root water uptake is considered as a volumetric sink term in the Richards equation \citep{molz1970extraction,FeddesR.A1978Sofw}. Macroscopic models have been favored in several studies \citep{molz1981models,mathur1999modeling,de2012root} due to their simplicity, as they do not necessitate detailed information about the root system's geometry. For instance, Molz and Remson \cite{molz1970extraction} emphasized the challenge of modeling water transport in soil when considering microscopic models. The dynamic and complex geometry of the root system makes it impractical to measure accurately. Additionally, the water permeability of roots varies with their positions along the root, which further complicates the modeling process \citep{kramer1969plant,molz1981models}. To address these complexities, most extraction functions have been developed using a macroscopic approach rather than a microscopic one \citep{molz1981models,mathur1999modeling,de2012root}.

In the present study, we use the macroscopic model proposed by Feddes et al. \citep{FeddesR.A1978Sofw} to represent plant root water uptake. This model accounts for plant water stress, normalized root distribution, and transpiration potential. In addition, the van Genuchten model \citep{van1980closed} is employed to describe capillary pressure, and the van Genuchten-Mualem model \citep{mualem1976new} is used to represent the relative permeability. For the validation and effectiveness of the proposed numerical model, we employed simplified models, such as stepwise and exponential forms \citep{yuan2005analytical} for plant root water uptake and the Gardner model for capillary pressure. The validation is established through comparison with existing analytical solutions. 

The design of suitable numerical schemes for the Richards equation present significant challenges due to to its highly nonlinear nature. In addition, incorporating an implicit sink term for root water uptake introduces additional numerical complexity of the governing equation. Few numerical methods have been proposed to solve the Richards equation with a  sink term due to plant root water uptake. For instance, Wilderotter \cite{wilderotter2003adaptive} used an adaptive finite-element method for solving the Richards equation with plant root growth. Janz and Stonier \cite{janz1995modeling} applied finite difference method to solve Richards equation with the presence of simplified model for water absorption by roots. Machado et al. \cite{machado2022new} used finite volume method to solve the Richards’ equation with evapotranspiration.  Efficient and robust numerical schemes are still in demand for modeling water flow in unsaturated soils while considering the root water uptake.

In terms of temporal discretization techniques, it is crucial to address the temporal derivative appropriately to ensure accurate numerical solutions for Richards' equation. This is particularly important due to the equation's highly nonlinear and the presence of stiff unsaturated soil properties \cite{berardi2016new}. Implicit or semi-implicit schemes are commonly employed for the temporal discretization to solve the Richards equation \citep{shahraiyni2012mathematical}. The primary reason for utilizing these schemes is to achieve stability and enable the use of practical time steps. The Backward Euler method is frequently used for temporal discretization in solving the Richards equation \citep{zha2019review}. 
Few second-order temporal schemes have been proposed in the literature for solving Richards equation, including the Crank-Nicolson method \citep{zha2016comparison,chavez2018numerical,cai2020convergence} and the backward differentiation formula \citep{baron2017adaptive,cumming2011mass,keita2021implicit}. 

The primary linearization methods applied to Richards' equation include the Newton method, Picard method, and the modified Picard method. The Newton method, which is quadratically convergent, has been successfully applied to Richards’ equation \cite{list2016study}. However, a significant drawback of the Newton method is its local convergence nature and the requirement for computing derivatives, which can be computationally intensive \cite{list2016study}.
The Picard technique, despite its wide usage in addressing Richards’ equation, faces challenges with convergence as highlighted in \cite{celia1990general,lehmann1998comparison}. An improvement is proposed in \cite{celia1990general} to obtain the modified Picard method. While this method maintains only linear convergence, it offers greater robustness compared to the Newton method \cite{list2016study}.

In this paper, an efficient approach that combines the localized exponential RBF method \citep{lee2003local} with a second-order backward differentiation formula for temporal discretization is developed. 
Localized RBF methods, which employ a set of scattered nodes distributed throughout the computational domain and its boundaries, eliminate the need for mesh generation and simplify the computational process. This local approach generates sparse matrices, which avoids ill-conditioning problems and reduces computational time \citep{lee2003local,boujoudar2021modelling,boujoudar2023localized}. The selection of EXP-RBF is due to its high convergence rate, as demonstrated in numerous studies \citep{madych1992bounds,cheng2003exponential,cheng2012multiquadric}. Additionally, its positive definiteness has been validated in various system solutions \citep{cheng2003exponential,fasshauer2007meshfree,cheng2012multiquadric}. Furthermore, the BDF2 scheme ensures stability and a reasonable choice of time steps. Furthermore, it enhances computational efficiency in dealing with stiffness in the resolution of the Richards equation  \citep{illiano2021iterative}. The modified Picard method \citep{celia1990general} is used to linearize the system. This technique has been used in various studies to overcome mass conservation problems encountered in solving the Richards equation \citep{celia1990general}. A variety of numerical schemes, such as the finite difference method \citep{dogan2005saturated}, the finite element method \citep{list2016study}, and the finite volume method \citep{manzini2004mass}, have been utilized to solve the Richards equation using the modified Picard iteration. 

In {\color{blue}this} study, the developed approach is enhanced by the use of modified Picard iterations. This approach is employed to solve the Richards equation with a sink term due to plant root water uptake. The main goal is to provide accurate numerical results of soil moisture in the root zone and to investigate the influence of plant root water uptake on the distribution of soil moisture.

The paper is structured as follows: Section \ref{secc:2} presents the mathematical model of infiltration and plant root water uptake. In Section \ref{sec:p4:3}, we introduce the LRBF meshless method and linearization techniques proposed to solve the governing system. Section \ref{sec:p4:4} is dedicated to numerical experiments for modeling soil moisture distribution in the root zone. Finally, Section \ref{sec:p4:5} provides the concluding remarks.
\section{Unsaturated flow and root water uptake models}
\label{secc:2}

\subsection{Unsaturated flow model} 
\label{secc:2.1}
Infiltration of water in soils is commonly described at the continuum scale using Richards' equation~\citep{richards1931capillary}. This equation can be derived by combining Buckingham-Darcy's law and the mass conservation equation \citep{bear1987theory}. The Richards equation with a sink term due to plant root water uptake and subject to initial and boundary conditions can be written as follows:  
\begin{equation}\label{E1:p5}
    \begin{cases}
          \dfrac{\partial\theta(\psi)}{\partial t}-\nabla.\left[K_\text{s}(\boldsymbol{x})k_\text{r}(\psi)\nabla \psi\right]-\dfrac{\partial \left[K_\text{s}(\boldsymbol{x})k_\text{r}(\psi)\right]}{\partial z}= -s(\boldsymbol{x},\psi),  \text{ $\Omega\times (0,T) $ }, \\
          -K_\text{s}(\boldsymbol{x})k_\text{r}(\psi)\nabla( \psi+1).n_{\Omega}=Q^N, \text{ $\partial\Omega^N\times (0,T), $ } \\
          \psi=\psi_D, \text{ $\partial\Omega^D\times (0,T), $ } \\
          \psi_{t=0}=\psi^0, \text{ $\Omega\times \lbrace 0\rbrace, $ }   
    \end{cases}
\end{equation}
where $\psi$ $[L]$ is the pressure head, $\theta$ $ [L^{3}/L^{3}]$ is the volumetric soil water content, $k_\text{r}(\psi)$ $ [-]$ is the water relative permeability, $K_\text{s}(\boldsymbol{x})$ $ [L/T]$ is the saturated hydraulic conductivity,
$\boldsymbol{x}=(x, y, z)^{T}$ is the coordinate vector, $x$ $[L]$ and $y$ $[L]$ denote the horizontal dimensions and $z$ $[L]$  denotes the vertical dimension positive upward, $T>0$ is a fixed time, $\Omega $ is an open set of $\mathbb{R}^{d}$ ($d = 1,2,3$), $\partial\Omega=\partial\Omega^D\cup\partial\Omega^N$ is the border of $\Omega$, $\psi_{D}$ represents the prescribed pressure head associated with the Dirichlet boundary $\partial \Omega^D$. Similarly, $Q^N$ denotes the infiltration flux associated to the Neumann boundary $\partial\Omega^N$, which can represent precipitation, irrigation, and evaporation rates, $\psi^0$ is the initial condition associated to the pressure head and $n_{\Omega}$ is the outward unit normal to the domain. The sink term $s( \boldsymbol{x},\psi)$ $ [L^{3}/L^{3} T]$ is defined as the volume of water removed per unit of time from a unit of volume of soil as a result of plant water uptake.  

Here, we use the mixed form of the Richards equation \eqref{E1:p5}, where both the pressure head and water content are considered. 
The Richards equation is highly nonlinear because of the nonlinear constitutive relations of the relative permeability and capillary pressure
in terms of saturation \citep{gardner1958some,brooks1964hydrau,van1980closed}.

The hydraulic conductivity $K=K_sk_r$ and water content in Equation \eqref{E1:p5} are estimated using empirical models, such as Gardner \citep{gardner1958some}, Brooks-Corey \citep{brooks1964hydrau}, van Genuchten \citep{van1980closed} and van Genuchten-Mualem \citep{mualem1976new,van1980closed}  models. In this study, the capillary pressure is described using the van Genuchten constitutive relationship \citep{van1980closed} model:
\begin{equation}\label{E2}
    \theta(\psi)=
    \begin{cases}
            \theta_r+\dfrac{\theta_s-\theta_r}{\left[1+(\alpha_{vg}\vert \psi \vert)^{n_{vg}}\right]^{m_{vg}}}, & \text{if}~\psi< 0, \\
            \theta_s, & \text{if}~\psi\geq 0,
    \end{cases}
\end{equation}
and the van Genuchten-Mualem model \citep{mualem1976new,van1980closed} is used for relative permeability:
\begin{equation}
K(\psi)=
\begin{cases}
    K_s \dfrac{\left[1-(\alpha_{vg}\vert \psi \vert)^{n_{vg}-1}\left[1+(\alpha_{vg}\vert \psi \vert)^{n_{vg}}\right]^{-m_{vg}}\right]^2}{\left[1+(\alpha_{vg}\vert \psi \vert)^{n_{vg}}\right]^{m_{vg }/2}}, & \text{if}~\psi< 0, \\
            K_s, & \text{if}~\psi\geq 0,
\end{cases}\label{E3:p5}
\end{equation}
 where $\theta_r$ and $\theta_s$ represent the residual and saturated water contents, respectively $[L^3/L^3]$, $\alpha_{vg}$ $[L^{-1}]$ is related to the inverse of the air-entry pressure, $n_{vg}>1$ is a measure of the pore-size distribution and $m_{vg}$ is given by $1-1/n_{vg}$.
 \subsection{Root water uptake model} 
\begin{figure}[ht!]
\centering
\includegraphics[width=8cm,height=5.5cm,angle=0]{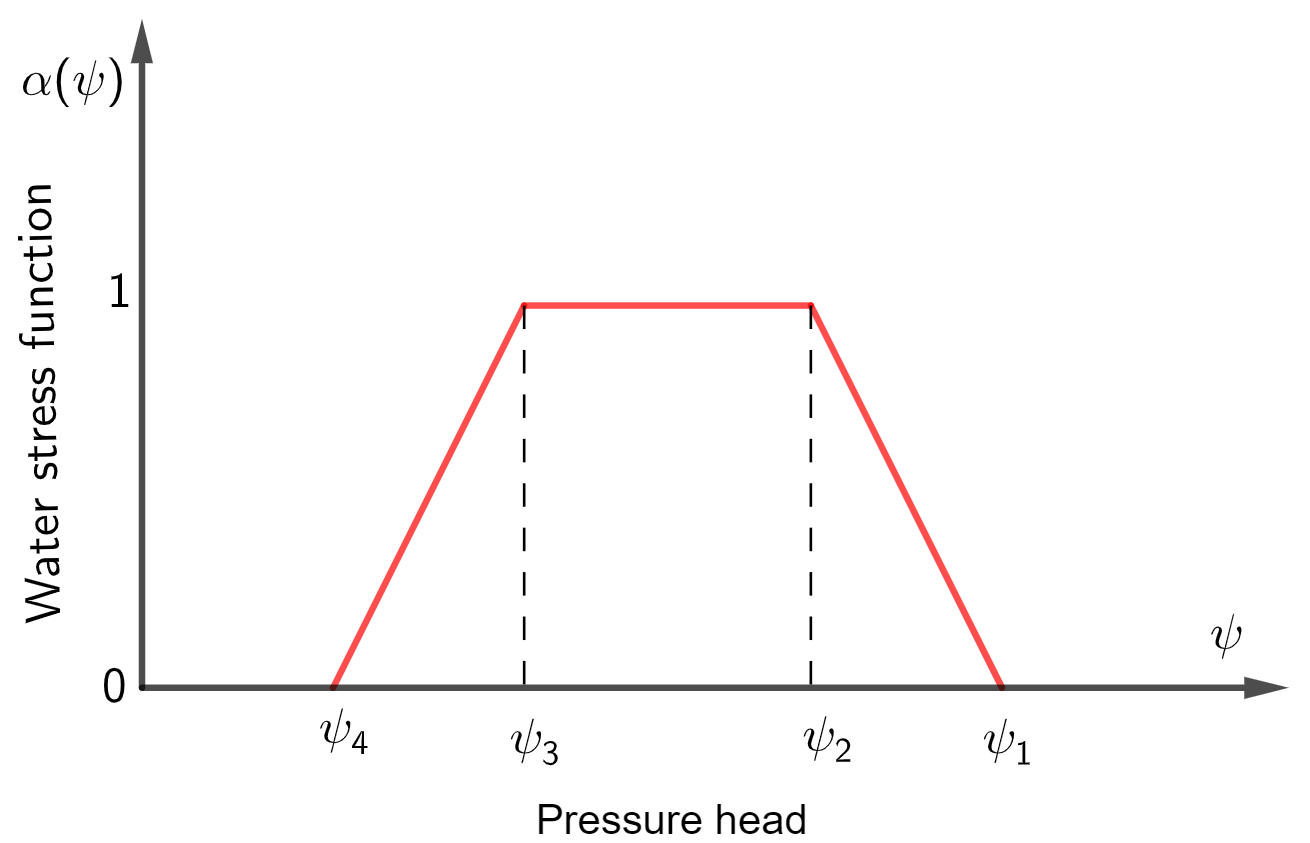} 
\caption{Schematic of the plant water stress response function \citep{FeddesR.A1978Sofw}.}\label{f1:p4}
\end{figure}

Various mathematical models are developed to describe plant root water uptake and they are based on microscopic or macroscopic approaches \citep{molz1981models}. In this study, we consider the macroscopic model proposed by Feddes et al. \cite{FeddesR.A1978Sofw}, which is expressed as follows:
\begin{equation}\label{E4:p5}
s( \boldsymbol{x},\psi)=\alpha(\psi)S_{p},
\end{equation}
where $ S_{p} $ is the potential water uptake rate $ [T^{-1}]$ and $\alpha(\psi) $ $ [-]$ is the soil water stress function ($ 0\leqslant \alpha(\psi) \leqslant 1$). Figure \ref{f1:p4} shows a schematic of the stress response function proposed by Feddes et al. \cite{FeddesR.A1978Sofw}, defined by the following equation:
\begin{equation}
\alpha(\psi)=
\begin{cases}
  0, &~\psi\geqslant \psi_{1}, \psi\leqslant \psi_{4}, \\
            \dfrac{\psi-\psi_{1}}{\psi_{2}-\psi_{1}}, &~\psi_{2}\leqslant \psi \leqslant \psi_{1}, \\
           1, &~\psi_{3}< \psi \leqslant \psi_{2}, \\
          \dfrac{\psi-\psi_{4}}{\psi_{3}-\psi_{4}}, &~\psi_{4}< \psi \leqslant \psi_{3}.
\end{cases}
\end{equation}

The stress response function $\alpha(\psi)$ is characterized by four critical pressure head values: $\psi_{1}$ denotes the pressure head level at which roots initiate water extraction from the soil, $\psi_{2}$ represents the pressure head level below which roots extract water at their maximum rate, $\psi_{3}$ indicates the limiting pressure head level below which roots cannot extract water at the maximum rate, and finally, $\psi_{4}$ signifies the pressure head level below which root water uptake ceases, often corresponding to the wilting point.

The variable $S_{p}$ denotes the water uptake rate under conditions of no water stress when $\alpha(\psi)=1$. The potential root water uptake $S_{p}$ is described according to Feddes et al.
 \cite{FeddesR.A1978Sofw} as follows:
\begin{equation}
S_{p}=b_r(\boldsymbol{x})T_{p},
\end{equation}
where $T_{p}$ denotes the potential transpiration rate $[LT^{-1}]$ and the function $b_r(\boldsymbol{x})$ $[L^{-3}]$ describes the spatial variation of the potential extraction term.
Various approaches are developed to express the spatial variation of the potential extraction term $b_r(\boldsymbol{x})$: constant, linear with the soil depth \citep{FeddesR.A1978Sofw}, or using the Hoffman and van Genuchten function \citep{hoffman1983soil}. For non-uniform root distributions, the function proposed by Vrugt \cite{vrugt2001one} can be employed.
\section{Numerical method}
\label{sec:p4:3}
In this study, an efficient approach is proposed to solve the governing equation.
This approach combines the localized EXP-RBF method with a second-order backward differentiation formula for temporal discretization. The modified Picards iteration method is employed to linearize the system. The localized EXP-RBF approach yields a sparse matrix for the global system, which reduces memory usage and computational time. The BDF2 scheme ensures the choice of reasonable time steps and enhances computational efficiency in dealing with stiffness in the resolution of the Richards equation. 
\subsection{Time discretization and linearization procedure}
Here, we employ a second-order backward differentiation formula to approximate equation \eqref{E1:p5}. The time interval $[0, T]$ is discretized into $N_t$ equally spaced intervals using a fixed time step size of $\Delta t = T/N_t$. Each discrete time point $t_n$ is defined as $t_n = n\Delta t$, where $n$ ranges from $0$ to $N_t$. The BDF2 discretization of equation \eqref{E1:p5} is defined as follows:
  \begin{equation}\label{E13:p5}
  \dfrac{3\theta^{n+1}-4\theta^{n}+\theta^{n-1}}{2\Delta t}-\nabla.(K^{n+1}\nabla \psi^{n+1})-\dfrac{\partial K^{n+1}}{\partial z}= -s^{n+1}.
  \end{equation}
The nonlinear equation \eqref{E13:p5} is solved by computing successive approximations of $\psi^{n+1}$. If $m$ denotes the iteration level, the Picard iteration scheme can be expressed as follows: 
    \begin{equation}\label{E14:p5}
  \dfrac{3\theta^{n+1,m+1}-4\theta^{n}+\theta^{n-1}}{2\Delta t}-\nabla.(K^{n+1,m}\nabla \psi^{n+1,m+1})-\dfrac{\partial (K^{n+1,m})}{\partial z}= -s^{n+1,m}.
  \end{equation}
In this study, we employ the mixed form of the Richards equation \eqref{E1:p5}, which incorporates both the $\psi$ and $\theta$ variables. We adopt the approach proposed by Celia et al. 
 \cite{celia1990general} to represent $\theta^{n+1,m+1}$ as a truncated Taylor series expansion with respect to $\psi$. The expansion can be written as:
  \begin{equation}\label{E15:p5}
    \theta^{n+1,m+1}=\theta^{n+1,m}+\dfrac{d\theta}{d \psi}\Big\vert^{n+1,m}\left(\psi^{n+1,m+1}-\psi^{n+1,m}\right)+0[(\delta^{m})^{2}],
  \end{equation}
    where $\delta^m$ is the iteration increment given by:
  \begin{equation}
      \delta^m=\psi^{n+1,m+1}-\psi^{n+1,m}.
  \end{equation}
Substituting Equation \eqref{E15:p5} into \eqref{E14:p5} yields the following expression:
    \begin{equation}\label{E17:p5}
 \left( \dfrac{3}{2\Delta t}C^{n+1,m}\right)\delta^m-\nabla.(K^{n+1,m}\nabla \psi^{n+1,m+1})=\dfrac{\partial K^{n+1,m}}{\partial z}-\dfrac{3\theta^{n+1,m}-4\theta^{n}+\theta^{n-1}}{2\Delta t}-s^{n+1,m},
  \end{equation}
  where $C=d\theta/d \psi$ is the specific moisture capacity function $[1/L]$. Equation \eqref{E17:p5} can be reformulated in terms of the iteration increment as follows:

\begin{align}\label{E18:p5}
\left( \dfrac{3}{2\Delta t}C^{n+1,m}\right)\delta^m-\nabla.(K^{n+1,m}\nabla \delta^{m}) &= \nabla.(K^{n+1,m}\nabla \psi^{n+1,m})+\dfrac{\partial K^{n+1,m}}{\partial z} \nonumber \\
&\quad -\dfrac{3\theta^{n+1,m}-4\theta^{n}+\theta^{n-1}}{2\Delta t}-s^{n+1,m},
\end{align}
 Equation \eqref{E18:p5} represents the general mixed-form Picard iteration, called the modified Picard approximation proposed by Celia et al. \cite{celia1990general}.  This representation is used in several studies for solving the Richards equation \citep{celia1990general,kirkland1992algorithms,huang1996new}. The final discrete form of the approximation can be obtained by applying finite difference method \citep{dogan2005saturated}, finite element method \citep{list2016study} and finite volume method \citep{manzini2004mass}. 
In this study, we use the localized EXP-RBF method to solve the Richards equation \eqref{E18:p5} based on the modified Picard iteration. 
\subsection{ Localized EXP-RBF method}
The localized EXP-RBF method is employed to solve the governing equation \eqref{E18:p5}. First, it is necessary to derive the two linearized operators: $\mathcal{L}^m$ corresponding to equation \eqref{E18:p5}, and $\mathcal{B}^m$ corresponding to the boundary conditions.

To simplify the three-dimensional representation of the spatial operator from Equation \eqref{E18:p5}, we introduce the following simplified expressions:
\begin{equation}
 \mathcal{L}_{d}^{m}\delta_{i}=\dfrac{\partial }{\partial x^{(d)}}\left( K_{i}\dfrac{\partial \delta_{i}}{\partial x^{(d)}} \right),   
\end{equation}
and
\begin{equation}
 \mathcal{L}_{4}^{m}K_{i}=\dfrac{\partial K_{i}}{\partial z},   
\end{equation}
where $(x^1,x^2,x^3)=(x,y,z)$ and $d$ represents the spatial dimension of the computation domain ($d = 1, 2, 3$).

To account for the heterogeneity of the porous medium and the dependence of $K$ on space, we propose to approximate the spatial operators $\mathcal{L}^{m}_{d}$ and $\mathcal{L}^{m}_{4}$ using the second-order standard finite difference method \citep{celia1990general}. This approach is based on the average ratio of $K$ at soil interfaces. The resulting approximations are presented in the following equations:
\begin{equation}
    \mathcal{L}_{d}^{m}\delta_{i}=\dfrac{1}{(\Delta x^{(d)})^{2}}\left(K^{(d)}_{i+1/2}(\delta^{(d)}_{iR}-\delta_{i})-K^{(d)}_{i-1/2}(\delta_{i}-\delta^{(d)}_{iL})\right).
\end{equation}
and
\begin{equation}\label{E22:p5}
   \mathcal{L}_{4}^{m}K_{i}=\dfrac{1}{\Delta z}\left(K^{(d)}_{i+1/2}-K^{(d)}_{i-1/2}\right),
\end{equation}
where $K^{(d)}_{i+1/2}$ and $K^{(d)}_{i-1/2}$ are given by:
\begin{equation}
    \begin{cases}
        K^{(d)}_{i+1/2}=\dfrac{1}{2}(K_{i}+K^{(d)}_{iR}), \\
         K^{(d)}_{i-1/2}=\dfrac{1}{2}(K_{i}+K^{(d)}_{iL}),
    \end{cases}
\end{equation}
where $K^{(d)}_{iR}$ and $K^{(d)}_{iL}$ represent the right and left neighboring points along the $x^{(d)}$-axis, respectively, which are given by:
\begin{equation}
    K^{(d)}_{iR}=\begin{cases}
        K(x_{i+1},y_i,z_i) & \text{if}~d=1, \\
         K(x_{i},y_{i+1},z_i) & \text{if}~d=2, \\
                K(x_{i},y_{i},z_{i+1}) & \text{if}~ d=3, \\
    \end{cases}
\end{equation}
and
\begin{equation}
    K^{(d)}_{iL}=\begin{cases}
        K(x_{i-1},y_i,z_i) & \text{if}~d=1, \\
         K(x_{i},y_{i-1},z_i) & \text{if}~d=2, \\
                K(x_{i},y_{i},z_{i-1}) & \text{if}~d=3. \\
    \end{cases}
\end{equation}
The same equations are utilized for the variable $\delta^m$ and $\psi^{n+1,m}$. The second spatial operator $\nabla \cdot (K_\text{s} k_\text{r} \nabla \delta^{m})$ in Equation \eqref{E18:p5} can be represented as follows:
\begin{equation}
\nabla.(K_\text{s}k_\text{r}\nabla \delta^{m})
=L^{m}\delta^m,
\end{equation}
where $L^{m}=(\mathcal{L}_{1}^{m}+\mathcal{L}_{2}^{m}+\mathcal{L}_{3}^{m}).$ The same applies for $\nabla.(K_\text{s}k_\text{r}\nabla \psi^{n+1,m})
=L^{m}\psi^{n+1,m}$.

Therefore, the Richards equation \eqref{E18:p5}, can be expressed in the following form:
  \begin{equation}
 \mathcal{L}^m \delta^{m}=a^{n+1,m}, 
  \end{equation}
  where:
  \begin{equation}
      \begin{cases}
      \mathcal{L}^m\delta^{m}=\left( \dfrac{3}{2\Delta t}C^{n+1,m}\right)\delta^m-L^m \delta^{m}, \\
      a^{n+1,m}=L^m \psi^{n+1,m}+\mathcal{L}_{4}^{m}K^{n+1,m}-\dfrac{3\theta^{n+1,m}-4\theta^{n}+\theta^{n-1}}{2\Delta t}-s^{n+1,m}.
      \end{cases}
  \end{equation}
$\mathcal{L}^m$ represents the linearized operator within the domain $\Omega$ at each time step $n+1$ and iteration $m$. 

Next, we will determine the linearized operator associated with the boundary conditions. Two boundary conditions, Dirichlet and Neumann, are considered and given by the following equations:
  \begin{equation}
      \begin{cases}
\psi^{n+1,m+1}=\psi_{D}, & \text{Dirchlet},\\
-K^{n+1,m}\nabla.(\psi^{n+1,m+1}+z).n_\Omega=Q^N, & \text{Neumann},
      \end{cases}
  \end{equation}
 These boundary conditions can be written in terms of the variable $\delta^m$. For Dirichlet conditions, it implies
that $\delta^m=0$. For Neumann conditions, it leads to $-K^{n+1,m}\nabla.(\delta^{m}+z).n_\Omega=Q^N-K^{n+1,m}\nabla.\psi^{n+1,m}.n_\Omega$. As a result, the boundary operator $\mathcal{B}^m$ can be represented as:
  \begin{equation}
\mathcal{B}^{m}\delta^{m}=b^{n+1,m},
  \end{equation}
  where
      \begin{equation}
      \mathcal{B}^{m}\delta^{m}=\begin{cases}
\delta^{m}, & \text{Dirchlet},\\
-K^{n+1,m}\nabla.(\delta^{m}+z).n_\Omega, & \text{Neumann}.
      \end{cases}
  \end{equation}
  and
  \begin{equation}
      b^{n+1,m}=\begin{cases}
0, & \text{Dirchlet},\\
Q^N-K^{n+1,m}\nabla.\psi^{n+1,m}.n_\Omega, & \text{Neumann}.
      \end{cases}
  \end{equation}
  Finally, the resulting system that needs to be solved is represented by:
    \begin{equation}\label{E32:p5}
      \begin{cases}
          \mathcal{L}^{m}\delta^{m}=a^{n+1,m},\\
      \mathcal{B}^{m}\delta^{m}=b^{n+1,m}.
      \end{cases}
  \end{equation}
 At each time level $t^{n+1}$, the linearized system \eqref{E32:p5} is solved using the localized EXP-RBF meshless method. The solution process entails iteratively solving the system at each iteration $m+1$ of the Picard iteration until the subsequent inequality is fulfilled:
  \begin{equation}\label{E33:p5}
      \vert \delta^{m} \vert \leq Tol,
  \end{equation}
where $Tol$ is the error tolerance.

\begin{figure}
  \centering
  \begin{tikzpicture}[scale=2,tdplot_main_coords]
    \draw[gray,fill=gray!20,opacity=0.7] (-1,-1,-1) -- (-1,1,-1) -- (1,1,-1) -- (1,-1,-1) -- cycle;
    \draw[gray,fill=gray!20,opacity=0.7] (-1,-1,1) -- (-1,1,1) -- (1,1,1) -- (1,-1,1) -- cycle;
    \draw[gray,fill=gray!20,opacity=0.7] (-1,-1,-1) -- (-1,-1,1) -- (1,-1,1) -- (1,-1,-1) -- cycle;
    \draw[gray,fill=gray!20,opacity=0.7] (-1,1,-1) -- (-1,1,1) -- (1,1,1) -- (1,1,-1) -- cycle;
    \draw[gray,fill=gray!20,opacity=0.7] (-1,-1,-1) -- (-1,-1,1) -- (-1,1,1) -- (-1,1,-1) -- cycle;
    \draw[gray,fill=gray!20,opacity=0.7] (1,-1,-1) -- (1,1,-1) -- (1,1,1) -- (1,-1,1) -- cycle;
    \coordinate (center) at (0,0,0);
    \draw[thick,fill=red] (center) circle (0.03);
    \node[below, red] at (center) {$\boldsymbol{x}_{s}$};
    \foreach \i/\label in {1/$\boldsymbol{x}_1^{[s]}$, 2/$\boldsymbol{x}_2^{[s]}$, 3/$\boldsymbol{x}_3^{[s]}$, 4/$\boldsymbol{x}_4^{[s]}$, 5/$\boldsymbol{x}_5^{[s]}$, 6/$\boldsymbol{x}_6^{[s]}$} {
      \pgfmathsetmacro{\theta}{\i*60}
      \coordinate (p\i) at ({cos(\theta)},{sin(\theta)},0);
      \draw[thick,fill=blue] (p\i) circle (0.03);
      \draw[dashed] (center) -- (p\i);
      \node[above, blue] at (p\i) {\label};
    } 
  \end{tikzpicture}
  \caption{Example of influence domain $\Omega_s$ with  seven points within a cube. }
  \label{f2:p4}
\end{figure}
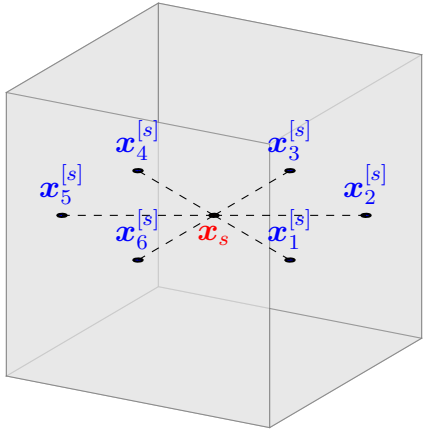

In the following, we describe the concept of the localized RBF method \citep{vsarler2006meshfree,li2013localized}. Let $\lbrace\boldsymbol{x_s}\rbrace_{s=1}^{N}$ represent a set of collocation points distributed in the domain $\Omega \cup\partial\Omega$, where $N$ denotes the total number of points. For each point $\boldsymbol{x}_s$ $\in\Omega \cup\partial\Omega$, we create a corresponding subdomain $\Omega_s=\left\lbrace \boldsymbol{x}_k^{[s]}\right\rbrace_{k=1}^{n_s}$ called influence domain, which contains the $n_s$ nearest neighbouring points to $\boldsymbol{x}_s$. We use the kd-tree algorithm \citep{bentley1975multidimensional} to determine the nearest neighbour points for each point $\boldsymbol{x}_s$. Figure \ref{f2:p4} displays an example of influence domain $\Omega_s$ with seven points within a cube.
The kd-tree algorithm is accurate and computationally efficient in identifying neighbouring nodes of the evaluation point \citep{yao2015implicit,li2013localized,boujoudar2021localized}. To approximate the variable $\delta^m$ using the localized EXP-RBF collocation method, we restrict the collocation scheme to the subdomain $\Omega_s$ instead of the entire domain $\Omega\cup\partial\Omega$. The solution $\delta^m$ can be approximated as a linear combination of RBFs within each influence domain as follows:
  \begin{equation}\label{E34:p5}
\delta^m(\boldsymbol{x}_s)=\sum_{k=1}^{n_{s}}\lambda^{[s]^{n+1,m+1}}_{k}\varphi(\Vert \boldsymbol{x}_{s}-\boldsymbol{x}_{k}^{[s]} \Vert),
 \end{equation}
 where $\{\lambda_{k}^{[s]^{n+1,m+1}}\}_{k=1}^{n_s}$ represent the unknown coefficients, and $\varphi$ denotes the radial basis function. In this study, we use the exponential RBF, given by:
\begin{equation}
\varphi(r_k)=\exp(-\varepsilon^2r_k^2),
\end{equation}
where $r_k=\Vert \boldsymbol{x}_{s}-\boldsymbol{x}_{k}^{[s]} \Vert$ denotes the distance between $\boldsymbol{x}_{s}$ and $\boldsymbol{x}_{k}^{[s]}$, and $\varepsilon>0$ is a shape parameter. The appropriate selection of the shape parameter plays a crucial role in ensuring the accuracy and stability of RBF meshless methods \citep{fasshauer2007choosing}. Various studies have proposed optimal choices for the shape parameter of specific RBFs \cite{hardy1971multiquadric,franke1998solving}.
The localized RBF methods are less influenced by the selection of the shape parameter compared to the global RBF methods, as shown in previous studies \citep{lee2003local,khoshfetrat2013numerical}.

The choice of the exponential RBF is motivated by its proven positive definiteness in various studies \citep{micchelli1984interpolation,cheng2012multiquadric,fasshauer2007meshfree}, which ensures that the resulting matrix is non-singular \citep{musavi1992training,garmanjani2018rbf}.  Additionally, a primary benefit of using the exponential RBF is the high convergence rate that can be achieved, as shown in several studies \citep{madych1992bounds,cheng2003exponential}.

Based on Equation \eqref{E34:p5}, we can represent the solution $\delta^{m}$ in matrix form as follows:
 \begin{equation}\label{E36}
\delta^{m}_{[s]}=\varphi^{[s]}\lambda^{[s]^{n+1,m+1}},
\end{equation}
where $\delta^m_{[s]}=\left[ \delta^m(\boldsymbol{x_{1}^{[s]}}), \delta^m(\boldsymbol{x_{2}^{[s]}}),..., \delta^m(\boldsymbol{x_{n_{s}}^{[s]}}) \right]^{T} $,\\ $\lambda^{[s]^{n+1,m+1}}=\left[ \lambda^{[s]^{n+1,m+1}}(\boldsymbol{x_{1}^{[s]}}), \lambda^{[s]^{n+1,m+1}}(\boldsymbol{x_{2}^{[s]}}),..., \lambda^{[s]^{n+1,m+1}}(\boldsymbol{x_{n_{s}}^{[s]}}) \right]^{T} $ and \\ $\varphi^{[s]}=\left[ \varphi(\Vert \boldsymbol{x}_{i}^{[s]}-\boldsymbol{x}_{j}^{[s]} \Vert ) \right]_{1\leq i,j\leq  n_s}$ is a real symmetric matrix of size $n_s \times n_s$. The unknown coefficients $\lambda^{[s]^{n+1,m+1}}$ can be determined as follows:
\begin{equation}\label{E37}
\lambda^{[s]^{n+1,m+1}}=(\varphi^{[s]})^{-1}\delta^m_{[s]}.
\end{equation}
Applying the linear operator $\mathcal{L}^{m}$ to Equation (\ref{E22:p5}) results in the following equations for $\boldsymbol{x_{s}} \in \Omega$:

\begin{equation} \label{E38}
\begin{split}
\mathcal{L}^{m}\delta^m(\boldsymbol{x_{s}}) & = \sum_{k=1}^{n_{s}}\lambda^{[s]^{n+1,m+1}}_{k}\mathcal{L}^{m}\varphi(\Vert \boldsymbol{x_{s}}-\boldsymbol{x_{k}}^{[s]} \Vert )=\sum_{k=1}^{n_{s}}\lambda^{[s]^{n+1,m+1}}_{k}\Psi^{m}(\Vert \boldsymbol{x_{s}}-\boldsymbol{x_{k}}^{[s]} \Vert ) \\ 
& = \vartheta_{[s]}^{m}\lambda^{n+1,m+1}_{[s]}=\vartheta_{[s]}^{m}(\varphi^{[s]})^{-1}\delta^m=\Lambda_{[s]}^{m}\mu_{[s]}^{n+1,m+1},
\end{split}
\end{equation}
where $\Psi^{m}=\mathcal{L}^{m}\varphi$, $\vartheta_{[s]}^{m}=\left[ \Psi^m(\Vert \boldsymbol{x_{s}}-\boldsymbol{x_{1}}^{[s]} \Vert),...,\Psi^m(\Vert \boldsymbol{x_{s}}-\boldsymbol{x_{n_{s}}^{[s]}} \Vert) \right]$ and $ \Lambda_{[s]}^{m}=\vartheta_{[s]}^{m}(\varphi^{[s]})^{-1} $. To reformulate Equation \eqref{E38} in terms of the global variable $\delta^{m}$ instead of the local variable $\delta_{[s]}^{m}$, we extend $\Lambda^{m}$ as the expansion of $\Lambda_{[s]}^{m}$ by adding zeros to the local vector where necessary. Thus, we establish the following correspondence:
\begin{equation} \label{E39}
    \mathcal{L}^{m}\delta_{[s]}^{m}(\boldsymbol{x_{s}}) =\Lambda^{m} \delta^{m},
\end{equation}
where $ \delta^{m}=\left[ \delta^{m}(\boldsymbol{x_{1}}),\delta^{m}(\boldsymbol{x_{2}}),...,\delta^{m}(\boldsymbol{x_{N}}) \right]^{T}$. 

Similarly, we apply the linear operator $\mathcal{B}^m$ to Equation \eqref{E34:p5} when $\boldsymbol{x_{s}}\in\partial\Omega$:
\begin{equation} \label{E40}
\begin{split}
\mathcal{B}^m\delta_{[s]}^{m}(\boldsymbol{x_{s}})&=\sum_{k=1}^{n_{s}}\lambda^{[s]^{m+1,n+1}}_{k}\mathcal{B}^m\varphi(\Vert \boldsymbol{x_{s}}-\boldsymbol{x_{k}}^{[s]} \Vert)=(\mathcal{B}^m\varphi^{[s]})\lambda^{[s]^{m+1,n+1}}\\
&=(\mathcal{B}\varphi^{[s]})(\varphi^{[s]})^{-1}\delta_{[s]}^{m}=\mathbb{\sigma}^{[s]}\delta_{[s]}^{m}=\mathbb{\sigma}\delta^{m}.
\end{split}
\end{equation}
Here, $\mathbb{\sigma}^{[s]}$ is defined as $(\mathcal{B}\varphi^{[s]})(\varphi^{[s]})^{-1}$, and $\mathbb{\sigma}$ represents the expansion of $\mathbb{\sigma}^{[s]}$ obtained by introducing zeros at the appropriate positions. By combining Equations \eqref{E39} and \eqref{E40} into Equation \eqref{E32:p5}, the following system is derived:
\begin{equation}\label{E41}
\begin{split}
\mathcal{L}^{m}\delta^{m}(\boldsymbol{x_{s}})=\Lambda^{m}(\boldsymbol{x_{s}})\delta^{m}=a^{n+1,m}(\boldsymbol{x_{s}}),\\
\mathcal{B}\delta^{m}(\boldsymbol{x_{s}})=\mathbb{\sigma}(\boldsymbol{x_{s}})\delta^{m}=b^{n+1,m}(\boldsymbol{x_{s}}).
\end{split}
\end{equation}
As a result, the localized EXP-RBF meshless method applied to the governing equation yields the final system to be solved:
\begin{equation}
    A^mX^m=B^m,
\end{equation}
where
$A^m=\left(
\begin{array}{c}
\Lambda^{m}(\boldsymbol{x_{1}}) \\
\Lambda^{m}(\boldsymbol{x_{2}}) \\
 . \\
.  \\
\Lambda^{m}(\boldsymbol{x_{N_{i}}})\\
\mathbb{\sigma}(\boldsymbol{x_{N_{i}+1}}) \\
 . \\
.  \\
\mathbb{\sigma}(\boldsymbol{x_{N}})
\end{array}\right)$, $X^m=\left(
\begin{array}{c}
\delta^{m}(\boldsymbol{x_{1}})\\
\delta^{m}(\boldsymbol{x_{2}})\\
.\\
.\\
\delta^{m}(\boldsymbol{x_{N_{i}}})\\
\delta^{m}(\boldsymbol{x_{N_{i}+1}})\\
. \\
.\\
\delta^{m}(\boldsymbol{x_{N}})
\end{array}
\right)$ and $B^m=\left(
\begin{array}{c}
a^{n+1,m}(\boldsymbol{x_{1}})\\
a^{n+1,m}(\boldsymbol{x_{2}})\\
.\\
.\\
a^{n+1,m}(\boldsymbol{x_{N_{i}}})\\
b^{n+1,m}(\boldsymbol{x_{N_{i}+1}})\\
. \\
.\\
b^{n+1,m}(\boldsymbol{x_{N}})
\end{array}
\right).$

This localized approach offers the advantage of inverting a sparse matrix, which mitigates the ill-conditioning issues encountered when dealing with the full matrix generated by the global method.

 By solving this sparse system, we obtain the approximate solution $\delta^{m}$ at the specified points. Once the condition specified in Equation \eqref{E33:p5} is verified, we assign $\psi^{n+1}$ the value of $\delta^{m}+\psi^{n,m+1}$. 
\section{Numerical experiments}\label{sec:p4:4}
In this section, we present numerical experiments for solving the Richards equation with sink terms due to plant root water uptake. The proposed approach uses the localized RBF meshless method based on the exponential function to solve the governing system \eqref{E32:p5}. To assess the effectiveness of the proposed numerical approach, several numerical tests are conducted. 

In the first test, two simplified models for root water uptake are used to validate the proposed numerical model by conducting a comparison between the numerical results and analytical solutions. In the second numerical test, the Feddes model \citep{feddes1976simulation} is employed to represent the plant root water uptake. The results obtained using the proposed numerical model are compared with those generated by Hydrus software \citep{simunek2005hydrus}, which serves as a benchmark for validation. The last numerical test evaluates the performance of the proposed numerical model to predict the three-dimensional soil moisture distribution profile in the root zone under axi-symetric irrigation conditions. The numerical tests are performed on a core CPU i7 2.1GHz computer in a MATLAB 2018a tools.

\subsection{Unsaturated flow in rooted soils under variable surface flux conditions}
In this numerical test, we validate our proposed numerical model by conducting a comparison between the approximate solutions and analytical solutions established in \citep{yuan2005analytical}. For this comparison, we use the Gardner model \citep{gardner1958some} to represent the relative permeability and a simple formulation for the water content \citep{warrick2003soil}, given by:
\begin{equation}
K=K_s \exp(\alpha \psi),
\end{equation}
\begin{equation}
\theta=\theta_r+(\theta_s-\theta_r) \exp(\alpha \psi).
\end{equation}
Both stepwise and exponential forms \citep{yuan2005analytical} are considered for root water uptake \citep{raats1974steady,rubin1993stochastic,yuan2005analytical}, which are respectively given by: 
\begin{equation}\label{E44'}
    s(z)=R_0\delta(z-l_1)=\begin{cases}
    R_0, & \text{if}~ l_1\leq z \leq l, \\
    0, & \text{if} ~0\leq z < l_1. \\
    \end{cases}
\end{equation}
\begin{equation}\label{E44}
    s(z)=R_0\exp\left[\beta(z-L)\right],
\end{equation}
where $\alpha$ $[1/L]$ is an empirical parameter related to the macroscopic capillary length of the soil, $R_0$ $[T^{-1}]$ is the maximum root uptake at the soil surface, $\beta$ $[L^{-1}]$ is a parameter represents the rate of reduction of root uptake, $l$ $[L]$ is the soil depth and $l_1$ $[L]$ is the maximum root depth. 

We perform numerical simulations of soil moisture and pressure head through a rooted soil with the following parameter values $l=100~cm$, $l_1=60~cm$, $\theta_s=0.45$, $\theta_r=0.2$, $K_s=1~cm~h^{-1}$ and $\beta=0.04~m^{-1}$ \citep{SrivastavaYeh,yuan2005analytical}. 

We consider two values of \( \alpha \) and \( R_0 \): \( \alpha = 0.01~\text{cm}^{-1} \) corresponds to \( R_0 = 0.02~\text{h}^{-1} \), and \( \alpha = 0.1~\text{cm}^{-1} \) corresponds to \( R_0 = 0.0025~\text{h}^{-1} \).
 The following boundary and initial conditions are used:
\begin{equation}\label{E45}
    \begin{cases}
    \psi(z,0)=\psi_0(z), \\
    \psi(0,t)=0, \\
   \left[K(\psi) \left( \dfrac{\partial \psi}{\partial z}+1 \right) \right ]_{z=L} =-q_1(t),
    \end{cases}
\end{equation}
where $\psi_0$ is the initial pressure head and $q_1$ is the flux at the soil surface, which changes over time. Our numerical test considers both steady-state and time-variant surface fluxes.

\begin{figure}[ht!]
\centering
\begin{tabular}{cc}\includegraphics[width=7cm,height=7.7cm,angle=0]{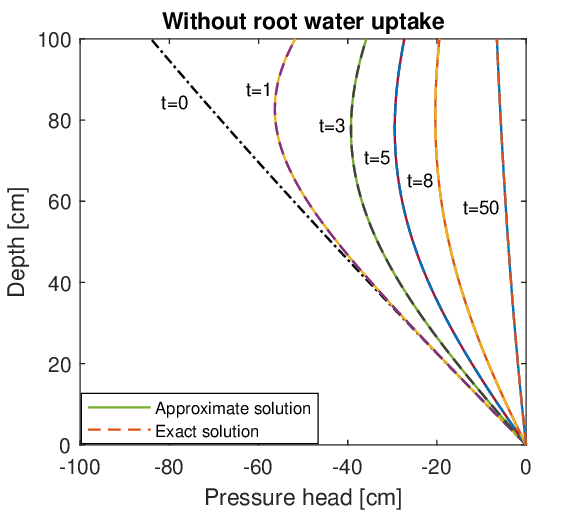}&
\includegraphics[width=7cm,height=7.7cm,angle=0]{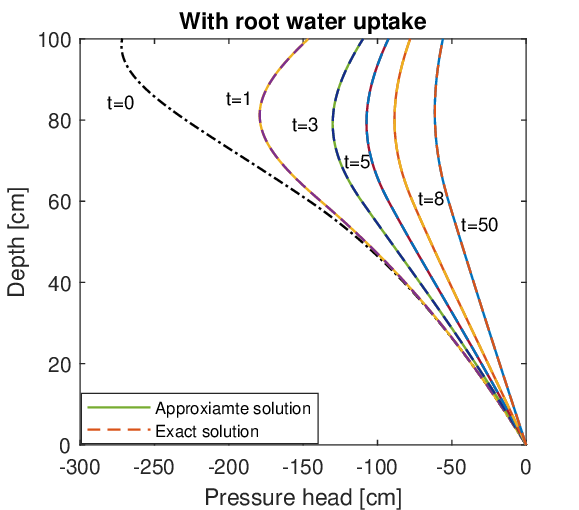}  \\
\includegraphics[width=7cm,height=7.7cm,angle=0]{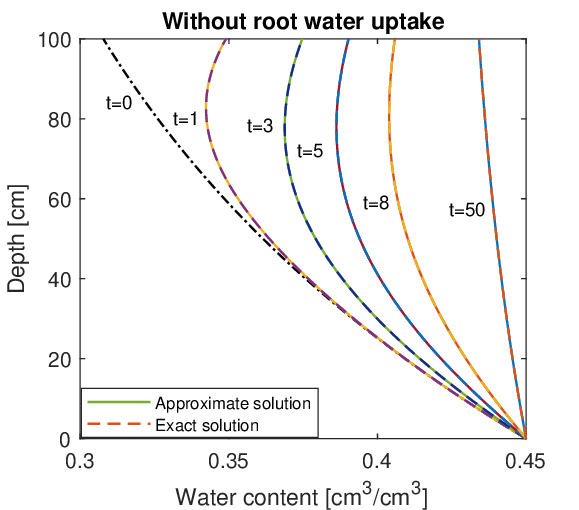} &\includegraphics[width=7cm,height=7.7cm,angle=0]{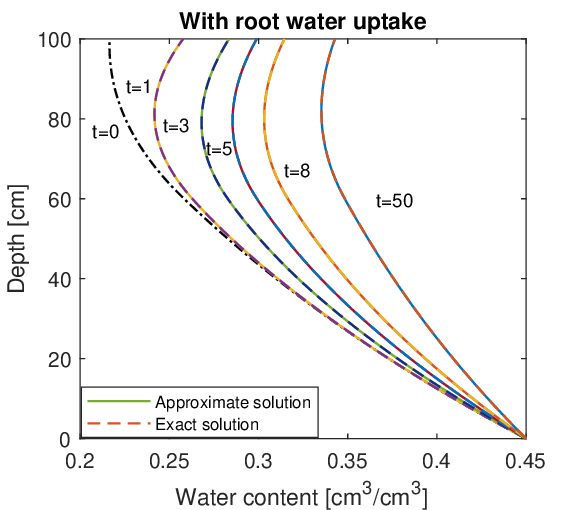}   
\end{tabular}
\caption{Case 1: $\alpha=0.01~cm^{-1}$ and $R_0=0.02~h^{-1}$. Comparison between approximate and exact solutions for water content and pressure head. Left: with root water uptake. Right: without root water uptake.}\label{f3:p5}
\end{figure}

\begin{figure}[ht!]
\centering
\begin{tabular}{cc}\includegraphics[width=7cm,height=7.7cm,angle=0]{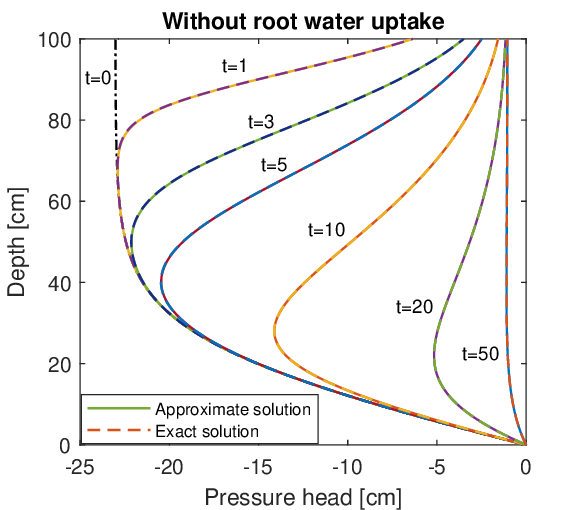}&
\includegraphics[width=7cm,height=7.7cm,angle=0]{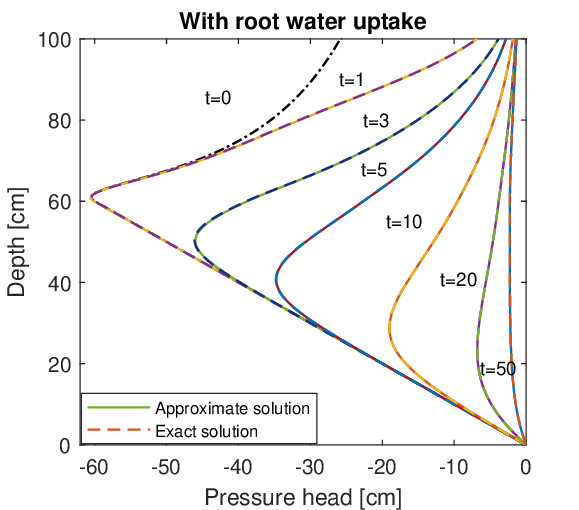}  \\
\includegraphics[width=7cm,height=7.7cm,angle=0]{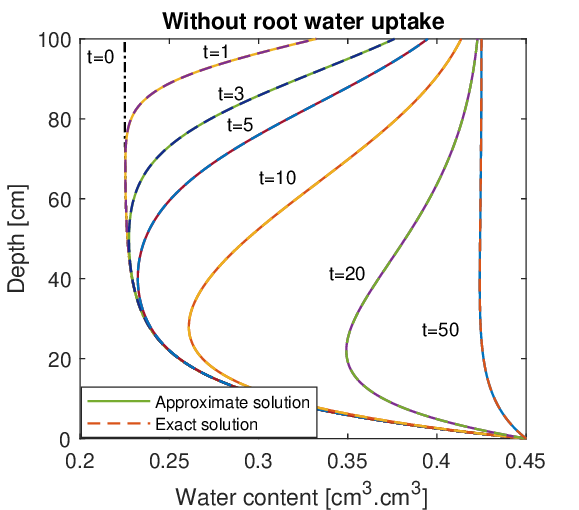} &\includegraphics[width=7cm,height=7.7cm,angle=0]{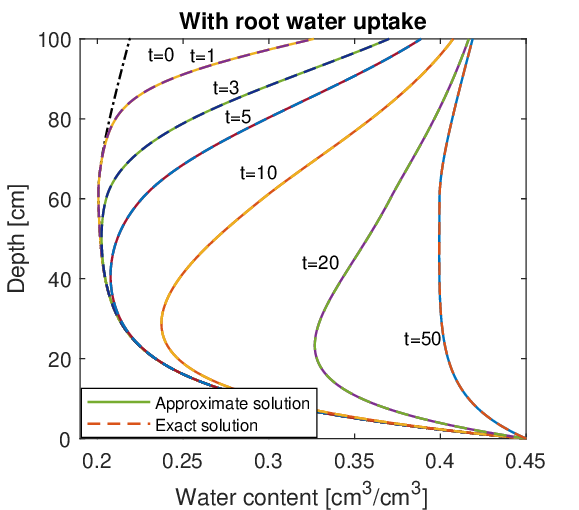}   
\end{tabular}
\caption{Case 1: $\alpha=0.01~cm^{-1}$ and $R_0=0.02~h^{-1}$. Comparison between approximate and exact solutions for water content and pressure head. Left: with root water uptake. Right: without root water uptake.}\label{f4:p5}
\end{figure}

First, we apply a steady-state surface flux as the upper boundary condition, considering a constant infiltration flux of $q_1=-0.9~cm~h^{-1}$. For the root water absorption, we use the formulation given by equation \eqref{E44'}. We present the distribution of water content and pressure head in Figures \ref{f3:p5} and \ref{f4:p5}, respectively. The simulations are performed using the two values of $\alpha$: $0.01~cm^{-1}$ and $0.1~cm^{-1}$.

The numerical simulations are conducted over a $50$-hour duration. To investigate the impact of root water uptake on soil moisture and pressure head distribution, we compare the time evolution of the water content and pressure head with and without considering root water uptake. The corresponding results are displayed for both cases, with and without root water uptake.

We observe a significant impact of root water uptake, particularly when $\alpha=0.01~m^{-1}$ and $R_0=0.02~h^{-1}$. The water content encounters significant changes under these conditions. However, for the case with $\alpha=0.1~m^{-1}$ and $R_0=0.0025~h^{-1}$, due to the relatively low value of the maximum water uptake, the effect on water content is relatively weak. 

\begin{figure}[ht!]
\centering
\begin{tabular}{cc}\includegraphics[width=7.4cm,height=6.4cm,angle=0]{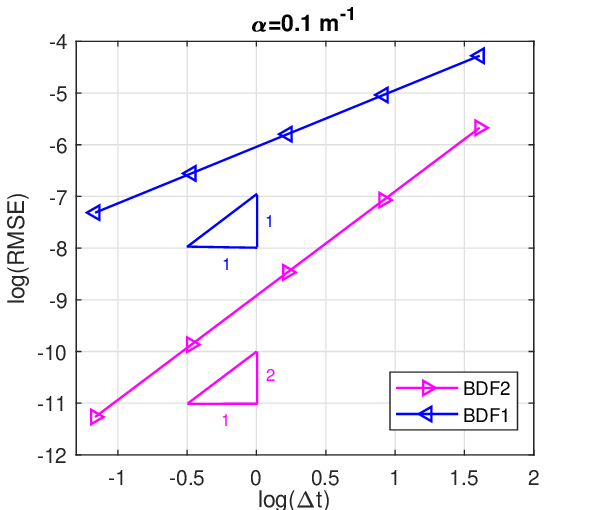}&
\includegraphics[width=7.4cm,height=6.4cm,angle=0]{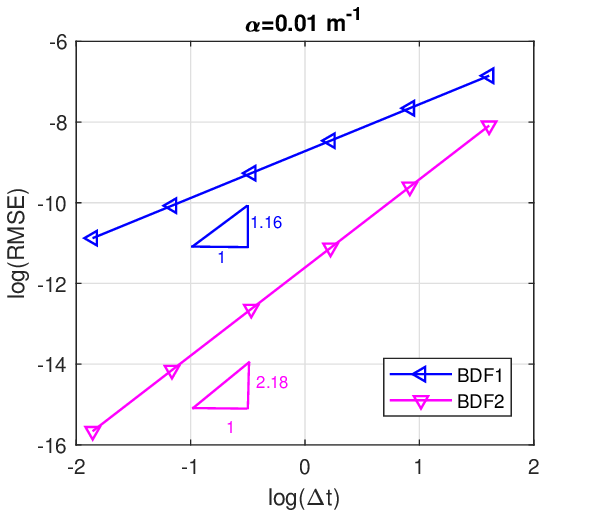}  \\
\end{tabular}
\caption{Root mean square error (RMSE) as a function of the time step for BDF1 and BDF2 schemes. Left: $\alpha=0.1~cm^{-1}$ and $R_0=0.0025~h^{-1}$. Right: $\alpha=0.01~cm^{-1}$ and $R_0=0.02~h^{-1}$. }\label{f5}
\end{figure}

With the same localized EXP-RBF parameters $N_z=1001$, $c=0.1$,  $\varepsilon=0.1$, and $\Delta t=0.01$, we display in Figure \ref{f5} the root mean square error (RMSE) associated to water content as a function of the time step for both first-order backward differentiation formula (BDF1) and BDF2 schemes.

 We confirm that the BDF2 scheme consistently achieves a convergence rate of approximately 2, while the BDF1 method exhibits a convergence rate of 1. Furthermore, BDF2 outperforms BDF1 in terms of accuracy.

\begin{figure}[ht!]
\centering
\begin{tabular}{cc}
\includegraphics[width=8cm,height=6cm,angle=0]{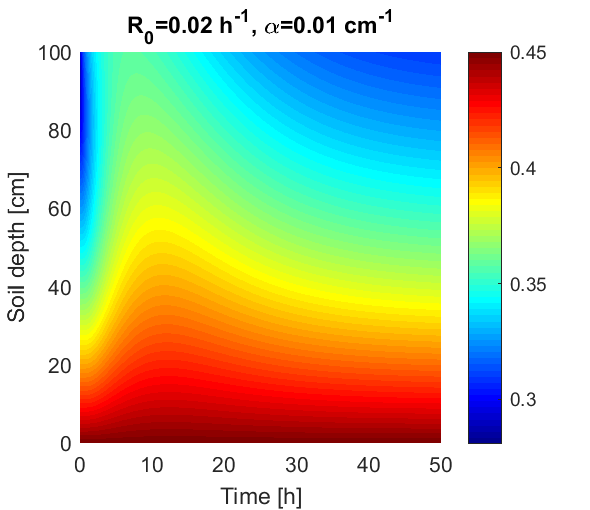} &
\includegraphics[width=8cm,height=6cm,angle=0]{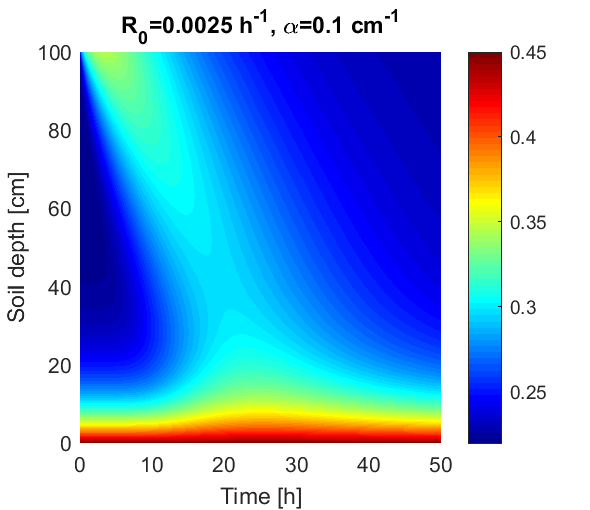}  
\end{tabular}
\caption{Time evolution of soil water content. Left: $\alpha=0.01~cm^{-1}$ and $R_0=0.02~h^{-1}$. Right: $\alpha=0.1~cm^{-1}$ and $R_0=0.0025~h^{-1}$.}\label{f6}
\end{figure}

Secondly, we introduce a time-varying surface flux as the upper boundary condition, providing a more realistic representation of real-world conditions influenced by factors such as evaporation, rainfall, and irrigation.
The flux at the upper boundary decreases exponentially with time, expressed as $q_1(t)=q_0+\delta \exp(k_1 t)$, with $\delta=-0.8~cm~h^{-1}$ and $k_1=-0.1~h^{-1}$. Root water uptake is modeled using the exponential formulation \eqref{E44}. Figure \ref{f6} illustrates the variation of water content in both time and space for the rooted soils.
The left figure corresponds to \( \alpha=0.01~\text{cm}^{-1} \) and \( R_0=0.02~\text{h}^{-1} \), while the right figure corresponds to \( \alpha=0.1~\text{cm}^{-1} \) and \( R_0=0.0025~\text{h}^{-1} \).

Both soils receive an equal water supply from the surface, but their moisture content patterns differ significantly. Soil profile 1, characterized by $R_0=0.02~h^{-1}$ and $\alpha=0.01~m^{-1}$, generally exhibits higher moisture levels than soil profile 2 with $R_0=0.0025~h^{-1}$ and $\alpha=0.1~m^{-1}$. Although soil profile 1 absorbs more water through its roots than profile soil 2, it benefits from the capillary rise, which facilitates water transfer from the water table to the root zones. This accounts for the noticeable difference in moisture content patterns between the two soils. 

Here, we will show the advantage of using BDF2, for temporal discretization in the proposed numerical model, in terms of computational time. To this end, we will analyze the CPU time of the proposed numerical techniques which are based on LRBF approach while using BDF2 or BDF1 for temporal integration. Numerical simulations are performed for the same numerical example until the final time $T=50$ hours for both temporal schemes. With a time step $\Delta t = 0.1$, BDF2 yields an RMSE of $1.64 \times 10^{-5}$ and requires a CPU time of $53$ seconds. For consistency in terms of comparison, we repeat the simulation using BDF1 in order to achieve the same level of accuracy which requires a time step $\Delta t = 0.015$ and leads to a CPU time of $293$ seconds. This comparison emphasizes the advantage of using BDF2 in our numerical model.

\begin{figure}[ht!]
\centering
\begin{tabular}{cc}
\includegraphics[width=5.7cm,height=5cm,angle=0]{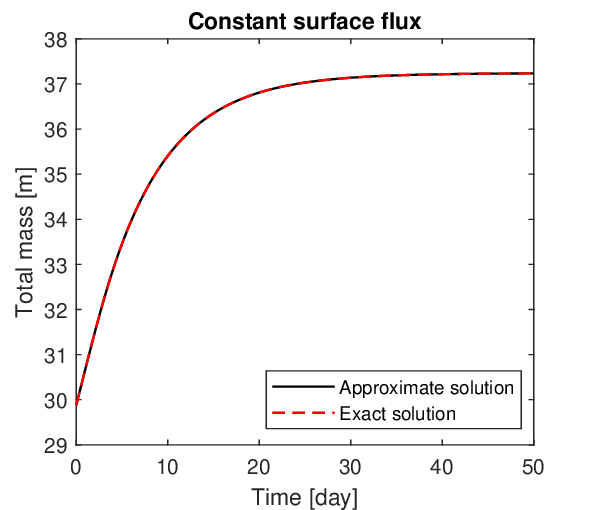}&
\includegraphics[width=5.7cm,height=5cm,angle=0]{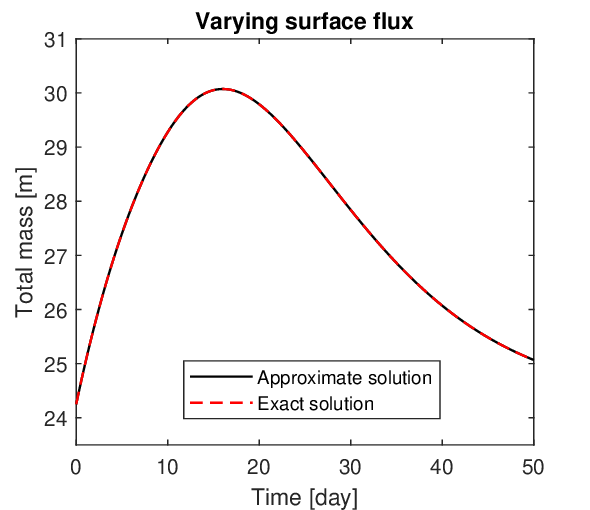}  
\end{tabular}
\caption{Comparison of the total mass between approximate and exact solutions. Left: Constant surface flux. Right: Varying surface flux.}\label{f7:rev}
\end{figure}

 The performance of the proposed numerical model in terms of total mass conservation is analyzed under both constant and varying surface flux conditions. We computed the evolution of the total mass of water ($I$) in the computational domain for the approximate and exact solutions.

\begin{equation}
    I(t)=\int_0^{L}\theta(z,t)dz,
\end{equation}
Figure \ref{f7:rev} shows the time evolution of total mass for the approximate and exact solutions. A good agreement is observed in both scenarios for constant and variable surface flux. This confirms the accuracy of the proposed numerical model in terms of mass conservation.

\subsection{Root water uptake soil
profile with groundwater table}
In this numerical example, we used Equations \eqref{E2},\eqref{E3:p5} and \eqref{E4:p5} to describe root water uptake. We consider a rooted soil with a depth of $1.2~m$ and a root depth of $0.9~m$. The root distribution is assumed to decrease linearly with depth. The bottom of the soil is in hydrostatic equilibrium with the groundwater table at the beginning of the simulation ($\psi(z,t=0)=-z$). The upper boundary condition is characterised by a non-flow boundary condition due to the absence of precipitation or irrigation. As a result, the plant's transpiration depletes the soil water, which is partially recharged by capillary action from the groundwater.

\begin{figure}[ht!]
\centering
\begin{tabular}{cc}
\includegraphics[width=7.3cm,height=7.3cm,angle=0]{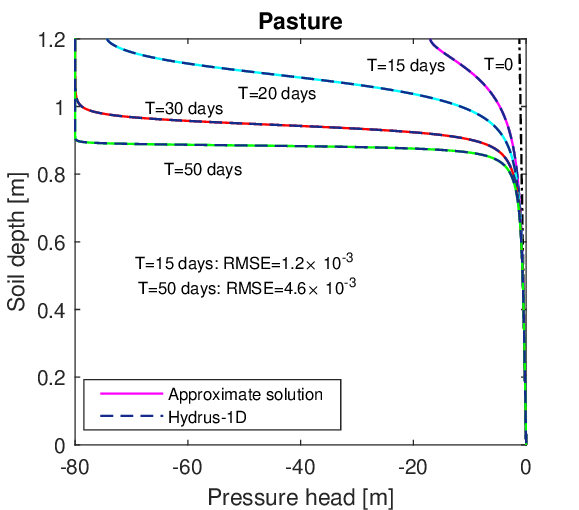}&
\includegraphics[width=7.3cm,height=7.3cm,angle=0]{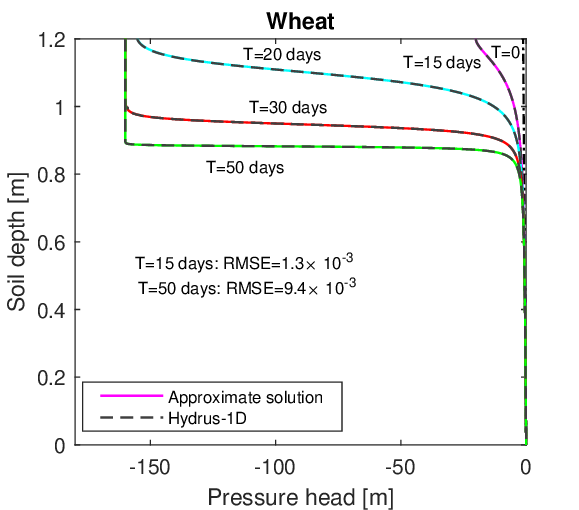}  
\end{tabular}
\caption{Pressure head profiles for two considered plants. Left: Pasture. Right: Wheat.}\label{f7}
\end{figure}

\begin{figure}[ht!]
\centering
\begin{tabular}{cc}
\includegraphics[width=7.3cm,height=7.3cm,angle=0]{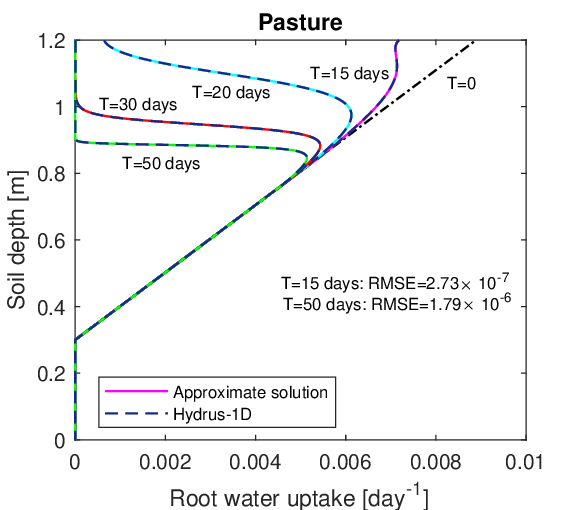}&
\includegraphics[width=7.3cm,height=7.3cm,angle=0]{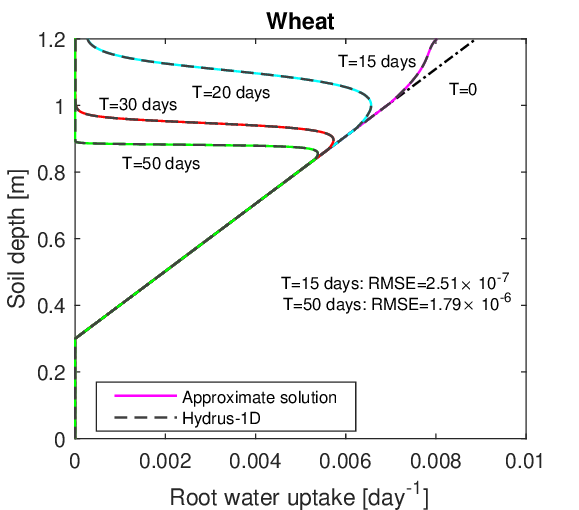}  
\end{tabular}
\caption{Root water uptake profiles for two considered plants. Left: Pasture. Right: Wheat.}\label{f8}
\end{figure}

Two types of plants, wheat and pasture, are considered, with their respective stress function parameters outlined in Table \ref{T1}. $\psi_{3,low}$ and $\psi_{3,high}$ represent the critical pressure head values, below which roots can no longer efficiently extract water at the maximum rates, $r_{2L}$ and $r_{2H}$, respectively. These potential transpiration rates, $r_{2L}$ and $r_{2H}$ are presently calibrated at $0.1~cm/day$ and $0.5~cm/day$, as reported in \citep{simunek2005hydrus}.
To get the value of $\psi_3$, we use the  following interpolation proposed by \cite{simunek2005hydrus}:
\begin{equation}\label{E46}
    \psi_3=\begin{cases}
    \psi_{3,\text{low}}+\dfrac{(h_{3,\text{high}}-h_{3,\text{low}})}{(r_{2L}-r_{2H})}(r_{2h}-T_p), ~&\text{$r_{2L}<T_p<r_{2H}$},\\
    \psi_{3,\text{low}}, ~&\text{$T_p\leq r_{2L}$},\\
     \psi_{3,\text{high}}, ~&\text{$T_p\geq r_{2H}$}.\\
    \end{cases}
\end{equation}
The hypothetical plants are presumed to exhibit a potential transpiration rate $T_{\text{p}}=4~mm/day$ \citep{SIMUNEK2009505,naghedifar2020numerical}.
\begin{table}[ht!]
\begin{center}
\caption{Plants parameters.}\label{T1}
 \begin{tabular}{|c|c|c|c|c|c|}
 \hline
 Plant & $\psi_1$ & $\psi_2$ & $\psi_{3,\text{low}}$& $\psi_{3,\text{high}}$ & $\psi_4$\\
  & $(m)$ & $(m)$ & $(m)$ & $(m)$ & $(m)$ \\
 \hline
 Pasture & $-0.1$& $-0.25$ &$-2$ & $-8$ & $-80$ \\ 
  \hline
 Wheat & $0$& $-0.01$ &$-5$ & $-9$ & $-160$ \\ 
  \hline
\end{tabular}
\end{center}
\end{table}
These simulations were conducted over a period of $50$ days. This numerical test is treated in \citep{SIMUNEK2009505,naghedifar2020numerical}, where the authors used Hydrus to simulate root water uptake. In this case, the van Genuchten and the van Genuchten–Mualem \citep{mualem1976new,van1980closed} models are used to describe the capillary pressure and the relative permeability, respectively. The soil investigated in this study is identified as loamy soil \citep{SIMUNEK2009505}, and its hydraulic properties are given by $\theta_s=0
.430$, $\theta_r=0.078$, $\alpha_{vg}=3.6~m^{-1}$, $K_s=0.2496~m/day$, $n_{vg}=1.56$ and $m_{vg}=0.3590$.
We selected $1001$ points along the $z$-axis and used a time step of $0.01$. The localized EXP-RBF function was configured with the parameters $\varepsilon=0.1$ and $n_s=3$. Figures \ref{f7} and \ref{f8} display the pressure head profiles and root water uptake at different time levels for two plants: pasture (left) and wheat (right), respectively. 
 A comparison is made between the approximate solutions obtained using the proposed numerical method and the results provided by Hydrus-1D. A good agreement is observed across all cases, indicating a good level of consistency. The accuracy of our proposed numerical method in predicting soil moisture, accounting for root water uptake, is reinforced by the obtained root mean square error (RMSE) values.
 
\begin{figure}[ht!]
\centering
\begin{tabular}{cc}
\includegraphics[width=7cm,height=6.5cm,angle=0]{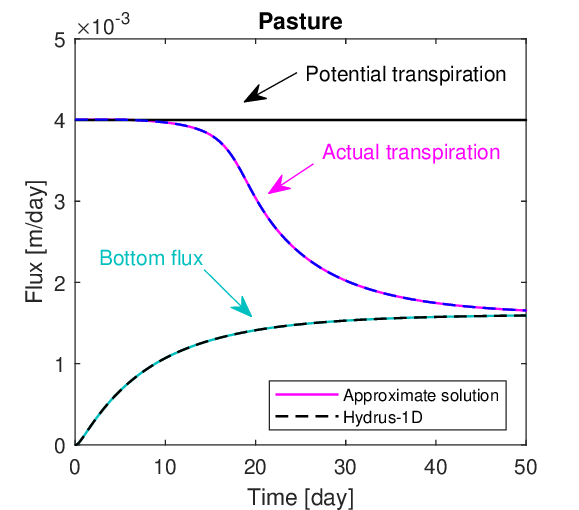}&
\includegraphics[width=7cm,height=6.5cm,angle=0]{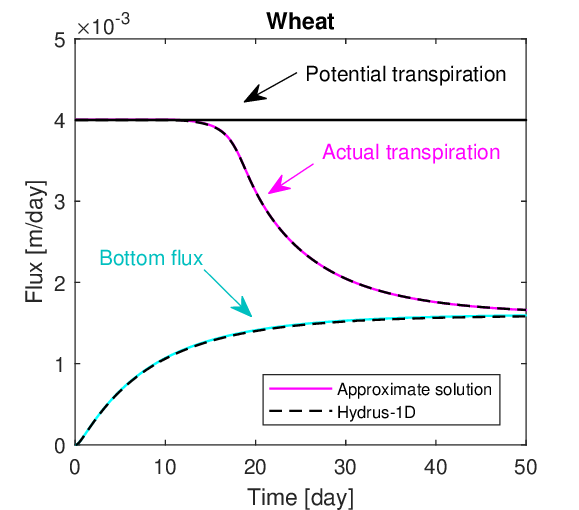}  
\end{tabular}
\caption{Time evolution of the potential transpiration, actual transpiration and bottom flux. Left: Pasture. Right: Wheat.}\label{f9}
\end{figure}

\begin{figure}[ht!]
\centering
\begin{tabular}{cc}
\includegraphics[width=7cm,height=6.5cm,angle=0]{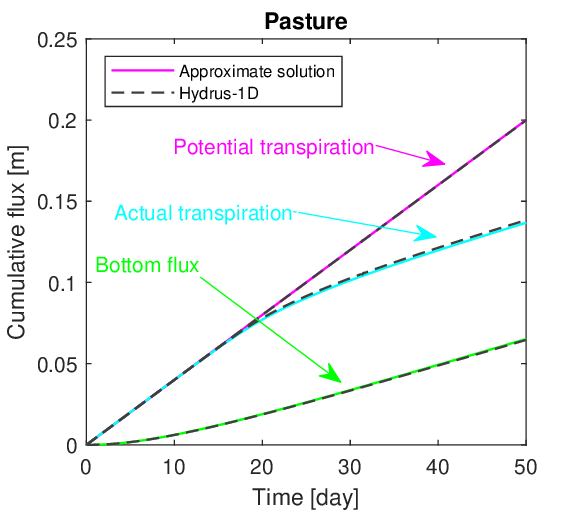}&
\includegraphics[width=7cm,height=6.5cm,angle=0]{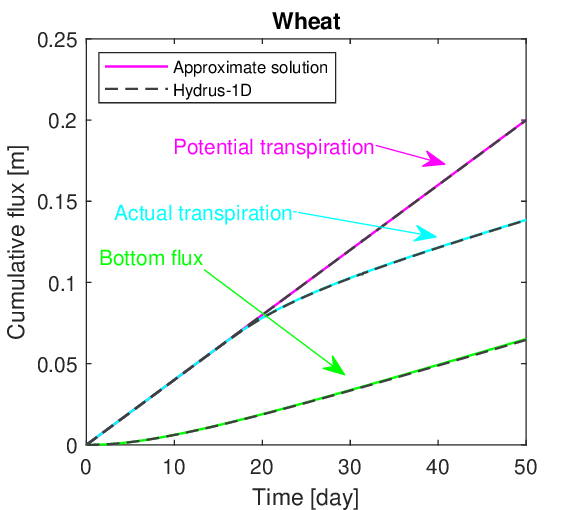}  
\end{tabular}
\caption{The cumulative flux of the potential transpiration, actual transpiration and bottom flux. Left: Pasture. Right: Wheat.}\label{f10}
\end{figure}

 The figures clearly demonstrate that RMSE values obtained from the numerical simulations are remarkably small, particularly when considering the error associated with the estimation of root water uptake. Figure \ref{f9}  presents the temporal evolution of the bottom flux $-K(\psi)\left[\partial{\psi}/\partial{z}+1 \right ]_{z=0}$, potential transpiration flux $T_p$, and actual transpiration flux $T_a$, as well as their respective cumulative values are displayed in Figure \ref{f10}. $T_a$ is given by \citep{simunek2005hydrus}:
\begin{equation}\label{E47}
    T_a=T_p \int_{\Omega_R} \alpha(\psi)b_r(z)dz.
\end{equation}

 These variations are observed over a period of $50$ days for the two plant types under investigation.

Once again, there is strong agreement between the approximate solutions using localized EXP-RBF meshless method and the results given by Hydrus-1D for all cases, demonstrating the effectiveness of the proposed numerical approach.

\begin{figure}[ht!]
\centering

\includegraphics[width=6.5cm,height=5.7cm,angle=0]{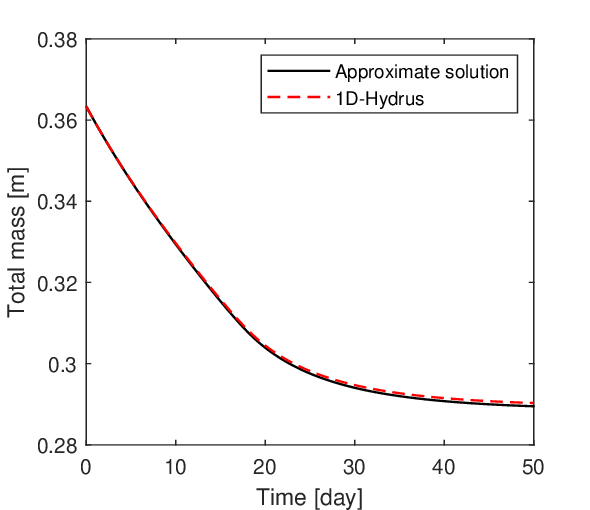}
\caption{Comparison of the total mass between the numerical and reference solutions.}\label{f11:rev}
\end{figure}

Figure \ref{f11:rev} presents the evolution of the total mass of water for the numerical solutions obtained using the proposed techniques and 1D-solution obtained using Hydrus. The results show the effectiveness of
the proposed numerical method in terms of conservation of mass.

\subsection{Three-dimensional soil moisture distribution profile under axi-symmetric irrigation conditions}
\begin{figure}[ht!]
\centering

\includegraphics[width=8cm,height=8cm,angle=0]{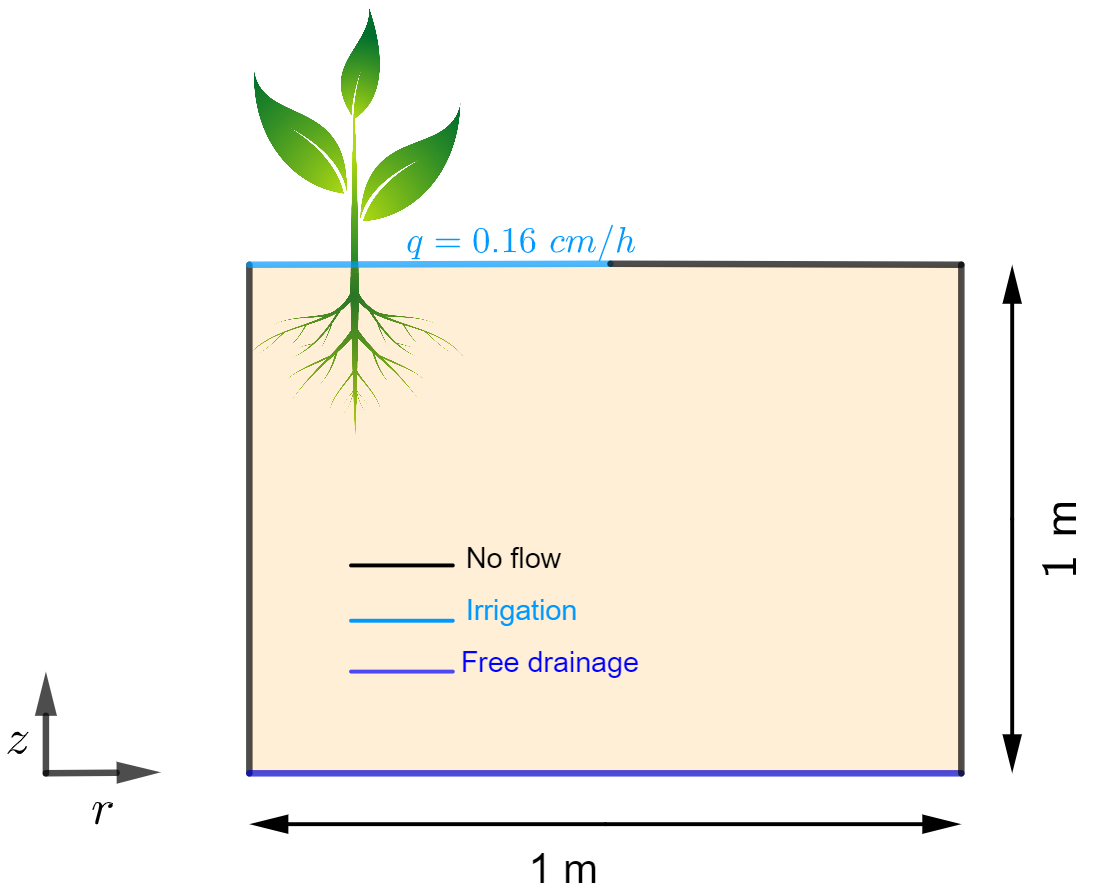}
\caption{Sketch of the computational domain.}\label{f11}
\end{figure}

This numerical test was used by Šimůnek and Hopmans \cite{SIMUNEK2009505}, where a three-dimensional axi-symmetrical profile is considered to simulate soil moisture in root zone. The computational domain can therefore be restricted to a square as shown in Figure \ref{f11}, where $r$ denotes the radial direction.
The domain's depth and radius are both set at $1~m$. An irrigation rate $q_z=0.16~cm/h$ is applied to the left half of the soil surface. In this case, the initial pressure head is set at $\psi_0=-4~m$ across the entirety of the soil profile, while the potential transpiration rate is set at $0.04~cm/h$. The Feddes model given by Equation \eqref{E4:p5} is used for the root water uptake model. The root density is assumed uniform in the radial direction and varied linearly with depth, from its maximum value at the soil surface to zeros at the soil depth $0.5 ~m$.

We conduct numerical simulations using the relationships described in equations \eqref{E2}–\eqref{E3:p5}, while employing the same parameters for the loamy soil as in the previous test. Furthermore, we use the same values for the soil water stress response parameters as listed in Table \ref{T1}. 
A zero flux boundary condition is imposed on the horizontal side of the computational domain to model the absence of flow across this boundary. Additionally, a free drainage condition is assumed at the bottom of the soil, representing the soil ability to freely drain water without any significant impedance or restriction. The boundary conditions are given by:
 \begin{equation}\label{E48}
   \begin{cases}
   -K \left(\dfrac{\partial \psi}{\partial z}+1\right)=-q_z, & z=1,\\
  -K\dfrac{\partial \psi}{\partial r}=0,
 & r=0, r=1, \\
   -K \dfrac{\partial \psi}{\partial z}=0, & z=0.\\
   \end{cases}
 \end{equation}
In \citep{SIMUNEK2009505}, the authors introduce two distinct sink terms to account for root water uptake. The first sink term corresponds to the conventional Feddes model, representing uncompensated root water uptake used in this study. The second sink term is known as compensated root water uptake. 

\begin{figure}[ht!]
\centering
\begin{tabular}{ccc}
\includegraphics[width=5.2cm,height=5cm,angle=0]{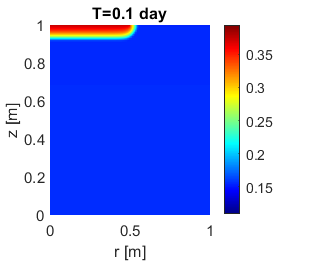} &
\includegraphics[width=5.2cm,height=5cm,angle=0]{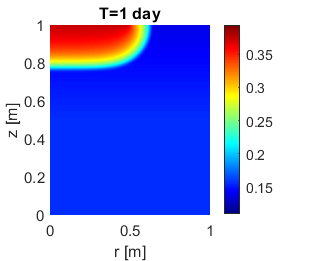} & \includegraphics[width=5.2cm,height=5cm,angle=0]{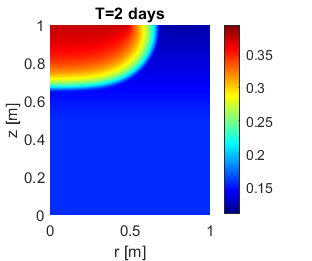}  \\
\includegraphics[width=5.2cm,height=5cm,angle=0]{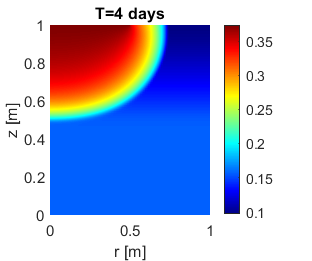} &
\includegraphics[width=5.2cm,height=5cm,angle=0]{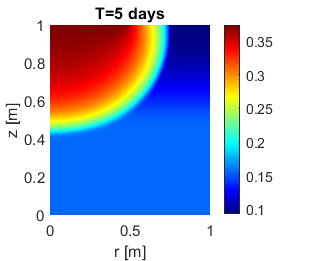} & \includegraphics[width=5.2cm,height=5cm,angle=0]{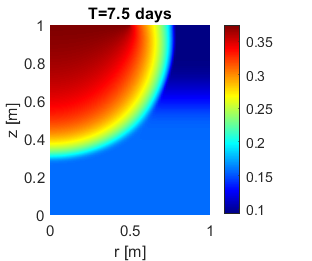}
\end{tabular}
\caption{The water content evolution at different times.}\label{f12:p4}
\end{figure}

In this numerical test, we use the following parameters: $\varepsilon=0.2$, $\Delta t=0.005$, $n_{s}=5$, $N_{x}=100$, and $N_{z}=200$. Figure \ref{f12:p4} illustrates the temporal evolution of water content in the case of uncompensated root water uptake. 

The numerical results align with the expected characteristics observed in the uncompensated scenario, as reported by \cite{SIMUNEK2009505}. 

\begin{figure}[ht!]
\centering
\includegraphics[width=8cm,height=7cm,angle=0]{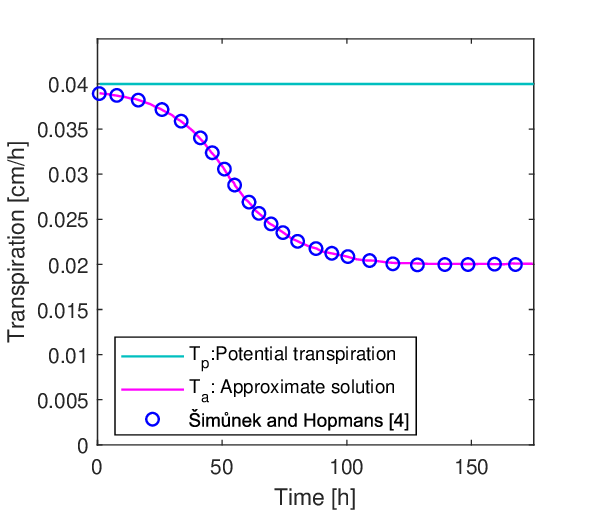}
\caption{Time evolution of the potential and actual transpiration.}\label{f13:p4}
\end{figure}

In Figure \ref{f13:p4}, we display the temporal evolution of the potential transpiration $T_p$ and actual transpiration $T_a$ computed using Equation \eqref{E47}. A comparison is performed between the approximate solutions obtained using the proposed numerical method and the results obtained from Hydrus \citep{SIMUNEK2009505}. The results demonstrate a good agreement, confirming the accuracy and reliability of our proposed numerical approach.

\section{Conclusion}
\label{sec:p4:5}
This study introduces efficient numerical techniques for solving the Richards equation with a sink term due to plant root water uptake. The van Genuchten and Feddes models are employed to describe capillary pressure and plant root water absorption, respectively. Our aim is to develop a numerical model that can accurately predict the soil moisture distribution in the root zone and to study the impact of plant root water absorption on soil moisture distribution.

An efficient approach is proposed, combining the localized exponential radial basis function method for space and the second-order backward differentiation formula for temporal discretization. The modified Picard's iteration method is used to linearize the system represented in the mixed form of the Richards equation.  The localized RBF methods, which use a unique set of scattered nodes distributed throughout the computational domain and its boundaries, eliminate the need for mesh generation and simplify the computational process. These techniques effectively overcome ill-conditioning problems often encountered in the use of global RBF methods since a sparse matrix for the global system is obtained. This localized approach improved memory usage and computational time. Additionally, the second-order accurate BDF2 scheme ensures reasonable time steps and enhances the accuracy of the method.

 The validation of the proposed numerical model for infiltration with plant root water absorption is performed using numerical experiments. The numerical results obtained using the proposed numerical model are compared against analytical and reference solutions.

The numerical model was validated for modeling unsaturated flow in rooted soils under variable surface flux conditions. Accurate results are obtained in comparison to the corresponding analytical solutions within a short CPU time. These results are confirmed by the RMSE values obtained across all considered cases, underlining the accuracy of the proposed numerical model in predicting the soil moisture dynamics in the root zone. Numerical simulations are performed using the proposed numerical techniques for modeling root water uptake soil
profile with groundwater table. We obtained good agreement between the numerical simulations and the results generated by the Hydrus software. The numerical simulations of three-dimensional soil moisture distribution profiles under axi-symmetric irrigation conditions, confirmed the accuracy and reliability of the proposed numerical model by comparing its findings with those presented in existing literature. The numerical simulations confirm that the combination of localized EXP-RBF and BDF2 methods leads to efficient numerical model for simulating soil moisture dynamics in the presence of root water uptake which is important in understanding soil-water-plant interactions.

\section*{\textbf{Declaration of competing interest}}
The authors declare that they have no known competing financial interests or personal relationships that could have appeared to
influence the work reported in this paper.

\section*{\textbf{Acknowledgments}}

AB gratefully acknowledges funding from UM6P-OCP. AB and AT gratefully acknowledge funding for the APRD research program from
the Moroccan Ministry of Higher Education, Scientific Research and
Innovation and the OCP Foundation.



\bibliographystyle{elsarticle-num}
\bibliography{mabiblio.bib}

\end{document}